# Vaccination strategies for SEIR models using feedback linearization. Preliminary results

M. De la Sen*, A. Ibeas** and S. Alonso-Quesada*

*Department of Electricity and Electronics. Faculty of Science and Technology
University of the Basque Country. PO. Box 644- Bilbao. Spain
** Departament of Telecommunications and Systems Engineering
School of Engineering. Universitat Autònoma de Barcelona. 08193 Bellaterra (Cerdanyola del Vallès). Spain

**Abstract**. A linearization-based feedback-control strategy for a SEIR epidemic model is discussed. The vaccination objective is the asymptotically tracking of the removed-by-immunity population to the total population while achieving simultaneously the remaining population (i.e. susceptible plus infected plus infectious) to asymptotically tend to zero. The disease controlpolicy is designed based on a feedback linearization technique which provides a general method to generate families of vaccination policies with sound technical background.

## 1. Introduction

Important control problems nowadays related to Life Sciences are the control of ecological models like, for instance, those of population evolution ( Beverton-Holt model, Hassell model, Ricker model etc.) via the online adjustment of the species environment carrying capacity, that of the population growth or that of the regulated harvesting quota as well as the disease propagation via vaccination control. In a set of papers, several variants and generalizations of the Beverton-Holt model (standard time–invariant, time-varying parameterized, generalized model or modified generalized model) have been investigated at the levels of stability, cycle- oscillatory behavior, permanence and control through the manipulation of the carrying capacity (see, for instance, [1-5]). The design of related control actions has been proved to be important in those papers at the levels, for instance, of aquaculture exploitation or plague fighting. On the other hand, the literature about epidemic mathematical models is exhaustive in many books and papers . A non-exhaustive list of references is given in this manuscript, cf. [6-14] (see also the references listed therein). Those epidemic models have also two major variants, namely, the so-called "pseudo-mass action models", where the total population is not taken into account as a relevant disease contagious factor and the so-called "true-mass action models", where the total population is more realistically considered as an inverse factor of the disease transmission rates). There are many variants of the above models, for instance, including vaccination of different kinds: constant [8], impulsive [12], discrete – time etc., incorporating point or distributed delays [12-13], oscillatory behaviours [14] etc. . On the other hand, variants of such models become considerably simpler for the illness transmission among plants [6-7]. In this paper, a feedback linearization technique is used to obtain a family of vaccination policies capable of asymptotically making the complete population become removed-by-immunity (immune). Initially, two of these vaccination policies are proposed and studied with detail. In a second stage, the general formalism is introduced to show the rationale behind the vaccination policies and how they fit into the general method. Feedback-linearization techniques, [17, 18], are successfully applied in control systems design for nonlinear problems, such as electrical machines, [19], or robotics [20]. However, its use in epidemic model control has been rather limited. It is assumed that the total population remains constant through time, so that the illness transmission is not critical, and the SEIR – model is of the above mentioned true-mass action type.

**1.1 Notation**. $R_+^n$ is the first open n-real orthant and $R_{0+}^n$ is the first closed n- real orthant.

$m \in R_{0+}^n$ is a positive real n-vector in the usual sense that all its components are nonnegative. In the same way, $M \in R_{0+}^{n \times n}$ is a positive real n-matrix in the usual sense that all its entries are nonnegative. The notations $R_+^n$ and $R_+^{n \times n}$ refer to the stronger properties that all the respective components or entries are positive.

$C^{(q)}(Do;Im)$ is the set of real functions of class q of domain Do and image Im. $PC^{(q)}(Do;Im)$ is the set of real functions of class (q-1) of domain Do and image Im whose q-th derivative exits but it is not necessarily everywhere continuous on its definition domain.

## 2. SEIR epidemic model

Let S (t) be the "susceptible" population of infection at time t, E (t) the " infected" ( i.e. those which incubate the illness but do not still have any symptoms) at time t, I (t ) is the " infectious" ( or "infective") population at time t, and R (t) is the " removed –by- immunity " ( or " immune") population at time t. Consider the SEIR-type epidemic model:

$$\dot{S}(t) = -\mu S(t) + \omega R(t) - \beta \frac{S(t)I(t)}{N} + \mu N(1 - V(t)) \qquad (1)$$

$$\dot{E}(t) = \beta \frac{S(t)I(t)}{N} - (\mu + \sigma)E(t) \qquad (2)$$

$$\dot{I}(t) = -(\mu + \gamma)I(t) + \sigma E(t) \qquad (3)$$

$$\dot{R}(t) = -(\mu + \omega)R(t) + \gamma I(t) + \mu N V(t) \qquad (4)$$

subject to initial conditions $S_0 = S(0) \geq 0$, $E_0 = E(0) \geq 0$, $I_0 = I(0) \geq 0$ and $R_0 = R(0) \geq 0$ under the vaccination constraint $V: R_{0+} \to R_{0+}$. In the above SEIR – model, N is the total population, μ is the rate of deaths from causes unrelated to the infection, ω is the rate of losing immunity, β is the transmission constant ( with the total number of infections per unity of time at time t being $\beta \frac{S(t)I(t)}{N}$), $\sigma^{-1}$ and $\gamma^{-1}$ are, respectively, the average durations of the latent and infective periods. All the above parameters are assumed to be nonnegative.

## 3. About the positivity of the SEIR epidemic model (1)-(4)

The vaccination strategy has to be implemented so that the SEIR model be positive in the usual sense that none of the populations, namely, susceptible, infected, infectious and immune be negative at any time. This requirement follows directly from the nature of the problem at hand. This section investigates conditions for positivity of the SEIR model (1)-(4) .The constant population constraint:

$$N = N(0) = S(t) + E(t) + I(t) + R(t) = S(0) + E(0) + I(0) + R(0) \ ; \ \forall t \in \mathbf{R}_{0+} \quad (5)$$

implying directly:

$$\dot{S}(t) + \dot{E}(t) + \dot{I}(t) + \dot{R}(t) = \dot{S}(0) + \dot{E}(0) + \dot{I}(0) + \dot{R}(0) = 0 \ ; \ \forall t \in \mathbf{R}_{0+} \quad (6)$$

holds directly in (1)-(4) since summing-up the four right- sides yields zero for all time.
The following assumption is made:

**Assumption 1**. The following constraints are assumed on the SEIR-model (1)-(4):
$min(S(0), I(0), R(0)) \geq 0$ and $E(0) > \frac{\mu+\gamma}{\sigma} I(0)$ with $\frac{\beta S(0) I(0)}{(\mu+\sigma)N} > E(0)$ if $I(0) \neq 0$. □

**Remark 1**. The physical interpretation of Assumption 1 is that the time origin of interest to fix initial conditions in the SEIR- model is the time instant at which the disease starts to be infectious. The growing rate of infectious at the time origin is positive, i.e. $\dot{I}(0) > 0$, even under zero initial condition $I(0) = 0$ so that $E(0) > 0$ and $\dot{E}(0) < 0$ for both the infected and infectious populations. Also, the infected population and its growing rate at the time origin are positive, i.e. $min(E(0), \dot{E}(0)) > 0$ if $I(0) = 0$ and a third consequence of Assumption 1 is that:

$\dot{E}(0) + \dot{I}(0) \leq 0$ if $\frac{I(0)}{E(0)} \leq \frac{\mu N}{\beta S(0) - (\mu+\gamma)N}$, $N \geq S(0) \geq \frac{\mu+\gamma}{\beta} N$ or, equivalently, if

$N \geq S(0) \geq \frac{\mu+\gamma}{N-\mu-\gamma}(E(0) + I(0) + R(0))$ requiring $\beta \geq \mu+\gamma$

$\dot{E}(0) + \dot{I}(0) > 0$ if $\frac{I(0)}{E(0)} > \frac{\mu N}{\beta S(0) - (\mu+\gamma)N}$, $N \geq S(0) \geq \frac{\mu+\gamma}{\beta} N$ requiring $\beta \geq \mu+\gamma$

so that $\dot{E}(0) + \dot{I}(0) \leq 0$ if $I(0) = 0$. Thus, if $I(0) = 0$ and $\dot{I}(0) > 0 (\dot{I}(0) \geq 0)$ then $\dot{E}(0) < 0 (\dot{E}(0) \leq 0)$ □

**Remark 2**. Note that Assumption 1 implies from (3) that:
- $\dot{I}(0) > 0$
- $S(0) < N - R(0) - \left(1 + \frac{\mu+\gamma}{\sigma}\right) I(0)$ and $S(0) > \frac{(\mu+\gamma)(\mu+\sigma)N}{\sigma \beta}$ if $I(0) \neq 0$
- $\beta > \beta_0 := (\mu+\gamma)\left(1 + \frac{\mu}{\sigma}\right)$ if $I(0) \neq 0$ (since $S(0) < N$). □

The parametrical condition $\beta > \beta_0$ is of interest even if $I(0) = 0$ in order to make the SEIR model parameters independent of any set of admissible initial conditions.

**Theorem 1**. Assume a vaccination function $V \in PC^{(0)}(\mathbf{R}_{0+}; [0,1])$ and that the initial conditions satisfy Assumption 1. Then, all the solutions of the SEIR – model (1)-(4) satisfy $S(t), E(t), I(t), R(t) \in [0, N]; \forall t \in \mathbf{R}_{0+}$.

*Proof*: The constant population constraint (5) is used in (1), (3)-(4) to eliminate the infected population E(t) leading to :

$$\dot{S}(t) = -(\mu + \alpha)S(t) + \omega R(t) + \left(\alpha - \beta \frac{I(t)}{N}\right) S(t) + \mu N (1 - V(t)) \tag{7}$$

$$\dot{I}(t) = -(\mu + \gamma + \sigma) I(t) + \sigma (N - S(t) - R(t)) \tag{8}$$

$$\dot{R}(t) = -(\mu + \omega) R(t) + \gamma I(t) + \mu N V(t) \tag{9}$$

for any given real constant $\alpha \geq (\beta/N) \sup_{t \geq 0}(I(t))$. Such a constraint is guaranteed with $\alpha \geq \alpha_0 := \beta$ if $0 \leq I(t) \leq N$ for all $t \geq 0$. It is possible to rewrite (7)-(9) in a compact form as a dynamic system of state $x(t) = (S(t), I(t), R(t))^T$, output $y(t) = S(t) + R(t)$ and whose input is appropriately related to the vaccination function as $u(t) = (1 - V(t), V(t))^T$. This leads to the following set of identities:

$$\dot{x}(t) = \overline{A}(\alpha) x(t) + \mu N \overline{E}_{13} u(t) + \left(\left(\alpha - \beta \frac{I(t)}{N}\right) E_1 x(t) + \sigma N e_2\right) \tag{10.a}$$

$$= A(\alpha) x(t) + \mu N \overline{E}_{13} u(t) + \left(\left[\left(\alpha - \beta \frac{I(t)}{N}\right) E_1 - \sigma E_{13}\right] x(t) + \sigma N e_2\right) \tag{10.b}$$

$$= A(\alpha) x(t) + \mu N \overline{E}_{13} u(t) + \left(\left(\alpha - \beta \frac{I(t)}{N}\right) E_1 x(t) + \sigma (N - y(t)) e_2\right) \tag{10.c}$$

$$= A(\alpha) x(t) + \mu N e_3 V(t) + \left(\left(\alpha - \beta \frac{I(t)}{N}\right) E_1 x(t) + \sigma (E(t) + I(t)) e_2 + \mu N e_1 (1 - V(t))\right) \tag{10.d}$$

$$y(t) = e_{13}^T x(t) \tag{11}$$

where $e_i$ is the i-th unit Euclidean column vector in $\mathbf{R}^3$ with its i-th component being equal to one and the other two components being zero, $e_{ij}$ having the i-th and j-th components being one and the remaining one being zero, so that $e_{13}^T = (1, 0, 1)$, and

$$\overline{A}(\alpha) := A(\alpha) - \sigma E_{13} \quad ; \quad A(\alpha) := \begin{bmatrix} -(\mu + \alpha) & 0 & \omega \\ 0 & -(\mu + \gamma + \sigma) & 0 \\ 0 & \gamma & -(\mu + \omega) \end{bmatrix} \tag{12}$$

$$E_{13} := \begin{bmatrix} 0^T \\ e_{13}^T \\ 0^T \end{bmatrix} = \begin{bmatrix} 0 & 0 & 0 \\ 1 & 0 & 1 \\ 0 & 0 & 0 \end{bmatrix} \; ; \; \overline{E}_{13} := [e_1, e_3] = \begin{bmatrix} 1 & 0 \\ 0 & 0 \\ 0 & 1 \end{bmatrix} \; ; \; E_1 := \begin{bmatrix} e_1^T \\ 0^T \\ 0^T \end{bmatrix} = \begin{bmatrix} 1 & 0 & 0 \\ 0 & 0 & 0 \\ 0 & 0 & 0 \end{bmatrix} \tag{13}$$

Note that: 1) $A(\alpha)$ is a Metzler matrix [4-5] for any given $\alpha \in \mathbf{R}_{0+}$ so that the $C_0$ – semigroup $\Phi \in L(\mathbf{R}^3, \mathbf{R}^3)$ of infinitesimal generator $A(\alpha)$ can also be represented as a fundamental positive real matrix function $\Psi := (\Phi x)(t) \in PC^{(1)}(\mathbf{R}_{0+}, \mathbf{R}_{0+}^{3\times 3})$ for $x \in Dom(\Phi) \subset \mathbf{R}^3$ of solutions of (7)-(9) as a result defined by $\Psi(\alpha, t) = e^{A(\alpha)t}$; $\forall t \geq 0$. In addition, since $\Psi(\alpha, t)$ is a fundamental matrix, it is nonsingular so that $(\Psi(\alpha, t)x) \in \mathbf{R}_+^3$; $\forall x \in \mathbf{R}_{0+}^3$

2) $min(\sigma, \mu) \geq 0$, $\alpha \geq \beta \frac{I(t)}{N}$, $V(t) \in [0,1]$; $\forall t \in \mathbf{R}_{0+}$; $e_1, e_2, e_3, e_{13} \in \mathbf{R}_{0+}^3$; $E_1 \in \mathbf{R}_{0+}^{3\times 3}$

3) From Assumption 1 (see also Remark 1), $E(0) + I(0) \geq 0$ and $\dot{E}(0) + \dot{I}(0) > 0$ if $\frac{I(0)}{E(0)} > \frac{\mu N}{\beta S(0) - (\mu + \gamma)N}$, $N \geq S(0) \geq \frac{\mu + \gamma}{\beta} N$ provided that $\beta \geq \mu + \gamma$. From continuity of any solution of (1)-(4), it exists $t_1 > 0$ such that $E(t) + I(t) > 0$; $\forall t \in (0, t_1)$. Also, $\dot{E}(0) + \dot{I}(0) \leq 0$ if $\frac{I(0)}{E(0)} \leq \frac{\mu N}{\beta S(0) - (\mu + \gamma)N}$ and $N \geq S(0) \geq \frac{\mu + \gamma}{\beta} N$ requiring $\beta \geq \mu + \gamma$. Thus, if $I(0) = 0$, $\dot{I}(0) > 0$ and $N \geq S(0) \geq \frac{\mu + \gamma}{\beta} N$ then $\dot{E}(0) + \dot{I}(0) \leq 0 \Rightarrow \dot{E}(0) < 0$. Again from continuity arguments, it exists $t_1 > 0$ such that $E(t) < E(0)$; $\forall t \in (0, t_1)$. Then, one has from (10.d) that for any admissible initial condition $x(0) = (S(0), I(0), R(0))^T$, the unique solution on $[0, t_1)$ of (7)-(9) is:

$$\mathbf{R}_+^3 \ni x(t) = e^{A(\alpha)t}\left(x(0) + \int_0^t e^{-A(\alpha)\tau} m(\tau) d\tau\right); \quad \forall t \in [0, t_1] \tag{14}$$

since $\mathbf{R}_+^3 \ni x(t) = e^{A(\alpha)t} x(0)$; $\forall x(0) \in \mathbf{R}_{0+}^3$, and since $x(t) \in \mathbf{R}_+^3$ on $[0, t_1)$ implies that

$$\mathbf{R}_+^3 \ni E(t) = e^{-(\mu + \sigma)t}\left[E(0) + \frac{\beta}{N} e_1^T \left(\int_0^t e^{(\mu + \sigma)\tau} x(\tau) x^T(\tau) d\tau\right) e_2\right]; \quad \forall t \in [0, t_1] \tag{15}$$

$$\mathbf{R}_{0+}^3 \ni m(t) := \mu N e_3 V(t) + \left(\left(\alpha - \beta \frac{I(t)}{N}\right) E_1 x(t) + \sigma(E(t) + I(t)) e_2 + \mu N e_1 (1 - V(t))\right); \quad \forall t \in [0, t_1]$$

$$= m_1(t) + \sigma\left(E(t) + e_2^T x(t)\right) e_2 \tag{16}$$

where $\mathbf{R}_{0+}^3 \ni m_1(t) := m(t) - \sigma\left(E(t) + e_2^T x(t)\right) e_2$ and $\left(\sigma\left(E(t) + e_2^T x(t)\right) e_2\right) \in \mathbf{R}_{0+}^3$. Since $x(t_1) \in \mathbf{R}_+^3$ and $e^{-(\mu + \sigma)t} E(0) \in \mathbf{R}_+$; $\forall t \in \mathbf{R}_{0+}$, it exists $t_2 > t_1$ such that $E(t) \in \mathbf{R}_+$, $m(t) \in \mathbf{R}_+^3$, $x(t) \in \mathbf{R}_+^3$ (so that $S(t), I(t), R(t) \in \mathbf{R}_{0+}$), $\forall t \in [0, t_2]$. The above properties extend to $t \in \mathbf{R}_{0+}$ from the structures of (14)-(16). Furthermore, $\left(\liminf_{t \to \infty} x(t)\right) \in \mathbf{R}_{0+}^3$ and $\left(\liminf_{t \to \infty} E(t)\right) \in \mathbf{R}_{0+}$. Those relations also imply from (5) that $max(S(t), E(t), I(t), R(t)) \leq N$; $\forall t \in \mathbf{R}_{0+}$. □

**Remark 3**. Note that the SEIR-model is not guaranteed to be positive according to Theorem 1 in the sense of [15-16] since Assumption 1 establishes constraints on the initial conditions. □

**Corollary 1**. Theorem 1 still holds if $V(t) \in \left[ 0, 1 + \left( \alpha - \beta \frac{I(t)}{N} \right) \frac{S(t)}{\mu N} \right]$; $t \in \mathbf{R}_{0+}$.

*Proof*: It follows from the proof of Theorem 1 since $m(t) \in \mathbf{R}_{0+}$; $\forall t \in \mathbf{R}_{0+}$ from (10.d) and (12)-(13) under this modified vaccination constraint. □

**4. Equilibrium points, stability, instability and immunity tracking of the whole population**

This section is concerned with vaccination designs so that stability or instability are guaranteed. It is also discussed how the vaccination might be synthesized so that the whole population is matched via vaccination strategies by the immune population so that the susceptible, infected and infectious are zeroed. An important point to deal with these issues is to ensure that all the partial populations (i.e. susceptible, infected, infectious and immune) are nonnegative for all time so that the boundndness of the whole population guarantees that of the individual ones. For this purpose, the positivity of the models is an important property to be guaranteed by the choice of the vaccination control in [0, 1]. Initially, the vaccination policies are introduced being its effects analyzed while in the following Section 5, the general frame (based-on feedback linearization) will be commented.

**Remark 4**. The equilibrium points in the vaccination-free case, which are discussed in Appendix A, are not suitable since one of them is concerned with the whole population being susceptible while the other is concerned with not all the population being asymptotically converging to the removed-by–immunity, in general. Therefore, a suitable vaccination strategy is necessary. The ideal vaccination mechanism objective is to reduce to zero the numbers of susceptible, infected and infectious independent of their initial numbers so that the total population becomes equal to the removed- by -immunity population after a certain time. After inspecting (1) and (4), it becomes obvious that the constraint $V : \mathbf{R}_{0+} \to \mathbf{R}_+$ is necessary to decrease the time variation of the susceptible and to increase simultaneously that of the removed by immunity. □

The following elementary result follows from the SEIR mathematical model (1)-(4)

**Assertion 1**. The SEIR model (1)-(4) fulfils the constant population through time constraint $N(t) := S(t) + E(t) + I(t) + R(t) = N(0) = N_0 = N > 0$ irrespective of the vaccination strategy.

*Proof*: It follows immediately by summing-up both sides of (1) to (4) what leads to:
$\dot{N}(t) = \dot{S}(t) + \dot{E}(t) + \dot{I}(t) + \dot{R}(t) = \mu(N(t) - S(t) - E(t) - I(t) - R(t)) = 0$; $\forall t \in \mathbf{R}_{0+}$
so that $N(0) = S(0) + E(0) + I(0) + R(0) = N_0 \Rightarrow N(t) = S(t) + E(t) + I(t) + R(t) = N_0 = N$ for all $t \geq 0$. □

**Remark 2**. Note that Assertion 1 proves that the constant population through time is independent of the vaccination strategy so that it is independent of the ideal vaccination objective constraint $V : \mathbf{R}_{0+} \to \mathbf{R}_+$ as a result. For instance, in a biological war, the objective would be to increase the numbers of the infected plus infectious population for all time. For that purpose, the appropriate vaccination strategy is negative. □

An auxiliary control function may be defined in several ways involving the vaccination function which is really the manipulated variable. For instance, one might define the infected/removed by immunity coupling term z(t) and control u(t) as follows:

$$z(t) = \omega R(t) - \sigma E(t) \tag{17}$$

$$u(t) = z(t) - \mu N V(t) \tag{18}$$

Note that any required control u(t) can be achieved using a vaccination strategy

$$V(t) = \frac{z(t) - u(t)}{\mu N} = \frac{\omega R(t) - \sigma E(t) - u(t)}{\mu N} \tag{19}$$

Then, one gets from (1)-(4) and (19):

$$\dot{E}(t) + \dot{I}(t) = -\mu(E(t) + I(t)) + \left(\beta \frac{S(t)}{N} - \gamma\right) I(t) \tag{20}$$

$$\dot{S}(t) + \dot{E}(t) = -\mu(S(t) + E(t)) + \mu N + u(t) \tag{21.a}$$

$$= \mu(I(t) + R(t)) + u(t) \tag{21.b}$$

$$\dot{I}(t) + \dot{R}(t) = -\mu(I(t) + R(t)) - u(t) = -\left(\dot{S}(t) + \dot{E}(t)\right) \tag{22}$$

$$\dot{S}(t) + \dot{R}(t) = -\mu(S(t) + R(t)) + \left(\gamma - \frac{\beta S(t)}{N}\right) I(t) + \mu N \tag{23.a}$$

$$= \mu(E(t) + I(t)) + \left(\gamma - \frac{\beta S(t)}{N}\right) I(t) \tag{23.b}$$

$$= \mu E(t) + \left(\mu + \gamma - \frac{\beta S(t)}{N}\right) I(t) = -\left(\dot{E}(t) + \dot{I}(t)\right) \tag{23.c}$$

$$\dot{S}(t) + \dot{E}(t) + \dot{I}(t) = -\mu(S(t) + E(t) + I(t)) + \omega R(t) - \gamma I(t) + \mu N (1 - V(t)) \tag{24.a}$$

$$= \mu R(t) + \sigma E(t) - \gamma I(t) + u(t) \tag{24.b}$$

$$\dot{R}(t) = -(\mu + \omega) R(t) + \gamma I(t) + \mu N V(t) \tag{25.a}$$

$$= -\mu R(t) - \sigma E(t) + \gamma I(t) - u(t) \tag{25.b}$$

$$= -\left(\dot{S}(t) + \dot{E}(t) + \dot{I}(t)\right) \tag{25.c}$$

**Remarks 3**. **(1)** Note from Eqs. (1) and (4) that a vaccination strategy $V(t) > 0$ applied on a time interval makes the susceptible population to decrease and the removed- by –immunity population to increase in a parallel fashion. From (2), the infected growing rate decreases as the susceptible numbers decrease.

**(2)** Eqs. 21-22 show that for a certain control associated with a vaccination strategy if the growing rate of joined susceptible plus infected population decreases then that of the infectious plus removed- by - immunity increases and conversely. □

The fact that the total population of the SEIR model (1)-(4) remains constant (Assertion 1) makes it both uncontrollable- to-the origin and unreachable. Such a constraint is atypical in most of control problems since the role of the vaccination is to decrease to zero the numbers of susceptible, infected and infectious to make the removed – by –immunity population to asymptotically converge to the total population.

**Assertion 2**. The SEIR model (1)-(4) is unreachable and uncontrollable to the origin via any vaccination strategy.

*Proof*: Proceed by contradiction. Fix any desired final state $x^* := (S(t^*), E(t^*), I(t^*), R(t^*))^T$ at arbitrary finite time $t = t^*$ fulfilling the constraint $S(t^*) + E(t^*) + I(t^*) + R(t^*) > N$. From Assertion 1, the population remains constant equal to N so that $x^*$ is unreachable at any time t*. Thus, the SEIR model is unreachable. It is always trivially uncontrollable-to-the origin for arbitrary initial conditions for the total population. □

A simple way of defining an useful control function is one with the goal of decreasing appropriately the numbers of susceptible while including the nonlinear term involving the product S(t)I(t) of susceptible and infectious in (1). The following result is concerned with this matter. A subsequent linear feedback vaccination strategy, being proportional to the susceptible for all time is discussed.

**Theorem 2**. The following properties hold:

**(i)** Assume that the feedback control and its associated vaccination strategy are generated as follows:

$$u(t) = -g\,S(t); \quad g \geq 0 \tag{26.a}$$

$$V(t) = \frac{1}{\mu N}\left(\omega R(t) + \left(g - \frac{\beta I(t)}{N}\right)S(t) + \mu N\right) \tag{26.b}$$

$\gamma \neq \sigma$, $g \neq \sigma$ and $g \neq \gamma$. Then the whole population becomes asymptotically removed by immunity at an exponential rate. Furthermore, $\exists \lim_{t \to \infty} V(t) := 1 + \frac{\omega}{\mu}$.

**(ii)** Assume that the feedback control and its associated vaccination strategy accordingly are generated as follows:

$$u(t) = -g\,(S(t) + E(t)), \quad g \geq 0 \tag{27.a}$$

$$V(t) = \frac{1}{\mu N}\left(g\,S(t) + (g - \sigma)E(t) + \omega R(t)\right) \tag{27.b}$$

$$= \frac{1}{\mu N}\left(g\,(N - I(t)) - \sigma E(t) + (\omega - g)R(t)\right) \tag{27.c}$$

Then $\lim_{t \to \infty} (S(t) + E(t)) = \frac{\mu N}{\mu + g}$ and $\lim_{t \to \infty} (I(t) + R(t)) = \frac{g N}{\mu + g}$ at an exponential rate if $0 \leq g < \mu$ and, furthermore,

$$\exists \lim_{t \to \infty}\left(V(t) + \frac{\sigma E(t) - \omega R(t)}{\mu N}\right) = \frac{g}{\mu + g} < 1$$

with $\lim_{t\to\infty} V(t) = \lim_{t\to\infty} E(t) = \lim_{t\to\infty} I(t) == \lim_{t\to\infty} R(t) = 0$ at exponential rates if $g = 0$. In particular, the whole population becomes asymptotically susceptible at an exponential rate if $0 = g < \mu$ for the corresponding vaccination law:

$$V(t) = \frac{1}{\mu N}\left(\omega R(t) - \sigma E(t)\right)$$

while all the other partial populations converge asymptotically to zero at an exponential rate.

*Proof*: **(i)** Rewrite (1) in the equivalent form:

$$\dot{S}(t) = -\mu S(t) + u(t) \qquad (28)$$

with an auxiliary control $u(t)$ being defined and generated as follows:

$$u(t) = \omega R(t) - \frac{\beta}{N} S(t) I(t) + \mu N(1 - V(t)) := -g S(t) \qquad (29)$$

through the vaccination function V(t) given by (26.b). Note that the open-loop solution of (28) is

$$S(t) = e^{-\mu t}\left(S(0) + \int_0^t e^{\mu \tau} u(\tau) d\tau\right) \qquad (30)$$

One gets from (28)-(29):

$$\dot{S}(t) = -(\mu + g) S(t) \Rightarrow S(t) = e^{-(\mu+g)t} S(0) \to S(\infty) = 0 \text{ as } t \to \infty \qquad (31)$$

$$u(t) = -g S(t) = -g e^{-(\mu+g)t} S(0) \to 0 \text{ as } t \to \infty \qquad (32)$$

Also, one gets from (31)-(32) into (29) that $\lim\sup_{t\to\infty}\left(V(t) - \frac{\omega}{\mu}\frac{R(t)}{N}\right) = 1$

implying that $V(t) \leq 1 + \frac{\omega}{\mu}\frac{R(t)}{N} + \varepsilon$; $\forall t \geq T = T(\varepsilon)$ (finite) and any arbitrary prefixed $\varepsilon \in \mathbf{R}_+$. On the other hand, one gets from (2) subject to (31):

$$\dot{E}(t) = -(\mu + \sigma) E(t) + \frac{\beta}{N} e^{-(\mu+g)t} S(0) I(t)$$

what leads to

$$E(t) = e^{-(\mu+\sigma)t} E(0) + \frac{\beta S(0)}{N} e^{-(\mu+\sigma)t} \int_0^t e^{-(g-\sigma)\tau} I(\tau) d\tau$$

$$\leq e^{-(\mu+\sigma)t} E(0) + \beta N \frac{e^{-(\mu+g)t} - e^{-(\mu+\sigma)t}}{\sigma - g} \to 0 \text{ exponentially fast as } t \to \infty \text{ at a rate of at}$$

most $\mu + \min(\sigma, g)$ as $t \to \infty$ since $\min(\mu, \sigma) > 0$ and $g \geq 0$. Combining this result with (3) and within the above solution expression for *E (t)*, one obtains:

$$I(t) = e^{-(\mu+\gamma)t} I(0) + \sigma \int_0^t e^{-(\mu+\gamma)(t-\tau)} E(\tau) d\tau$$

$$\leq e^{-(\mu+\gamma)t} I(0) + \sigma\, e^{-(\mu+\gamma)t} \int_0^t e^{(\mu+\gamma)\tau} \left( e^{-(\mu+\sigma)\tau} E(0) + \beta N \frac{e^{-(\mu+g)\tau} - e^{-(\mu+\sigma)\tau}}{\sigma - g} \right) d\tau$$

$$= e^{-(\mu+\gamma)t} I(0) + \frac{\sigma E(0)}{\gamma - \sigma}\left( e^{-(\mu+\sigma)t} - e^{-(\mu+\gamma)t} \right)$$

$$+ \frac{\sigma \beta N}{(\sigma-g)(\gamma-g)}\left( e^{-(\mu+g)t} - e^{-(\mu+\gamma)t} \right) - \frac{\sigma \beta N}{(\sigma-g)(\gamma-\sigma)}\left( e^{-(\mu+\sigma)t} - e^{-(\mu+\gamma)t} \right)$$

$$= e^{-(\mu+\gamma)t} I(0) + \frac{\sigma\,((\sigma-g)E(0) - \beta N)}{(\gamma-\sigma)(\sigma-g)} e^{-(\mu+\sigma)t}$$

$$+ \frac{\sigma\,(\beta N - (\gamma-g)E(0))}{(\gamma-\sigma)(\gamma-g)} e^{-(\mu+\gamma)t} + \frac{\sigma \beta N}{(\gamma-g)(\sigma-g)} e^{-(\mu+g)t} \to 0$$

exponentially fast as $t \to \infty$ provided that $\gamma \neq \sigma$, $g \neq \sigma$ and $g \neq \gamma$ and $g > \max(-\mu, -\sigma, -\gamma)$, which is guaranteed since $g \geq 0$, so that $R(t) = (N - (S(t) + E(t) + I(t))) \to N$ exponentially fast as $t \to \infty$.

Furthermore, from (26.b), $\lim\sup_{t \to \infty}\left( V(t) - \frac{\omega}{\mu} \frac{R(t)}{N} \right) = 1 \Rightarrow \exists \lim_{t \to \infty} V(t) = 1 + \frac{\omega}{\mu}$ since $R(t) \to N$ as $t \to \infty$. Property (i) has been proven.

**(ii)** Eqs. (15) together with (1)-(2) lead to:

$$\dot{S}(t) + \dot{E}(t) = -\mu(S(t) + E(t)) + \mu N + u(t) = -(\mu+g)(S(t) + E(t)) + \mu N \tag{33}$$

whose solution is subject to

$$S(t) + E(t) = e^{-(\mu+g)t}\left( S(0) + E(0) + \mu N \int_0^t e^{(\mu+g)\tau} d\tau \right) \to \frac{\mu N}{\mu + g} \text{ as } t \to \infty \tag{34}$$

satisfies $(I(t) + R(t)) \to \frac{g N}{\mu + g}$ as $t \to \infty$. As a result, $S(t) + E(t) \to N$ and $I(t), R(t) \to 0$ as $t \to \infty$ at exponential rate if $g = 0$, that is if the vaccination law is $V(t) = \frac{1}{\mu N}(\omega R(t) - \sigma E(t))$, since the SEIR – model (1)–(4) is a positive dynamic system if $g < \mu$ and $V(t) \in [0, 1]$; $\forall t \in \mathbf{R}_{0+}$ as it is proven in Appendix B. Then, all the partial populations are nonnegative for all time. From (2), $\exists \lim_{t \to \infty}\left( V(t) + \frac{\sigma E(t) - \omega R(t)}{\mu N} \right) := \frac{g}{\mu + g} < 1$. Also, $I(t) \to 0$ as $t \to \infty \Rightarrow E(t) \to 0$ as $t \to \infty$ for $g = 0$ so that $S(t) \to N$, since $I(t) \to 0$, $E(t) \to 0$ and $R(t) \to 0$ as $t \to \infty$ as $t \to \infty$ and $\lim_{t \to \infty} V(t) = 0$ from (27.b) or from (27.c). □

A more convenient vaccination strategy because of its properties and because of its implementation issues, based on measuring the immune and total populations instead of the susceptible one, is proposed in the subsequent result:

**Theorem 3**. Assume that the control and its associated vaccination strategy are as follows:

$$u(t) = -g R(t) + g_1 N, \quad g > -(\mu + \omega) \tag{35}$$

$$V(t) = \frac{1}{\mu N}(g_1 N - g R(t) - \gamma I(t)) \tag{36}$$

The above vaccination strategy implies that the removed- by- immunity population equalizes asymptotically the total population at exponential rate while the sum of the infected, infectious and susceptible populations converge asymptotically to zero at exponential decay rates according to:

$$R(\infty) = \frac{g_1 N}{\mu + \omega + g}; \quad S(\infty) + E(\infty) + I(\infty) = \frac{(\mu + \omega + g - g_1)N}{\mu + \omega + g} \tag{37}$$

$$\lim_{t \to \infty} \int_0^t e^{-\mu(t-\tau)}(\omega + g)R(\tau)d\tau = \left(\frac{g_1}{\mu} - \frac{g_1}{\mu + \omega + g}\right)N = \frac{g_1(\omega + g)}{\mu(\mu + \omega + g)}N \tag{38}$$

irrespective of the initial conditions and if, in particular, $g_1 = \mu + \omega + g$ then

$$R(\infty) = N; \quad S(\infty) + E(\infty) + I(\infty) = 0; \quad \lim_{t \to \infty} \int_0^t e^{-\mu(t-\tau)}(\omega + g)R(\tau)d\tau = \frac{(\omega + g)N}{\mu} = \frac{(g_1 - \mu)N}{\mu}$$

The vaccination effort is nonnegative for all time if $g_1 \geq \gamma$ or if $g \geq 0$ and $\gamma = \mu + \omega$ provided that the SEIR model supplies nonnegative populations.

*Proof*: Combining (23)-(24) yields a control action:

$$u(t) = \gamma I(t) + \mu N V(t) = -g R(t) + g_1 N \tag{39}$$

so that (4) becomes:

$$\dot{R}(t) = -(\mu + \omega)R(t) + u(t) = -(\mu + \omega + g)R(t) + g_1 N \tag{40}$$

so that if $g > -(\mu + \omega)$ then

$$R(t) \to R(\infty) := \lim_{t \to \infty} R(t) = \frac{g_1 N}{\mu + \omega + g} \quad \text{as } t \to \infty \tag{41}$$

at exponential rate according to an absolute upper-bound of exponential order of (41) being equal to $-(\mu + \omega + g) < 0$ and $S(\infty) + E(\infty) + I(\infty) = \frac{(\mu + \omega + g - g_1)N}{\mu + \omega + g}$. Also, $R(\infty) = N$ if $g_1 = \mu + \omega + g$ as a result. From Assertion 1, (41) with $g_1 = \mu + \omega + g$ implies also that $S(\infty) + E(\infty) + I(\infty) = 0$. On the other hand, the vaccination function is nonnegative from (36) for all time if $g_1 \geq \gamma$ since

$$g_1 N - g R - \gamma I = g_1(R + S + E + I) - g R - \gamma I = (g_1 - g)R + (g_1 - \gamma)I + g_1(S + E)$$
$$= (\mu + \omega)R + (g_1 - \gamma)I + g_1(S + E) \geq 0$$

If the constraint $g_1 \geq \gamma$ is changed to $g \geq 0$ and $\gamma = \mu + \omega$ then $g_1 - g = \gamma = \mu + \omega$, equivalently $g_1 - \gamma = g \geq 0$, and the vaccination function is also nonnegative for all time since the above expression becomes:

$$g_1 N - gR - \gamma I = (\mu + \omega)R + gI + g_1(S + E) \geq 0 \qquad \square$$

**Remark 4**. Note that Theorem 3 holds in particular if $g = 0$ in (35)-(36) so that the vaccination strategy is adopted on the basis of taking into account the total population only. Summing-up (1) and (2) and using (39) yields:

$$\dot{S}(t) + \dot{E}(t) = -\mu(S(t) + E(t)) + \omega R(t) - \sigma E(t) + \mu N(1 - V(t))$$
$$= -\mu(S(t) + E(t)) + (\omega + g)R(t) + \gamma I(t) - \sigma E(t) + (\mu - g_1)N \qquad (42)$$

which combined with (3) yields:

$$\dot{S}(t) + \dot{E}(t) + \dot{I}(t) = -\mu(S(t) + E(t)) + \omega R(t) - (\mu + \gamma)I(t) + \mu N(1 - V(t))$$
$$= -\mu(S(t) + E(t) + I(t)) + (\omega + g)R(t) + (\mu - g_1)N \qquad (43)$$

leading to the following solution of susceptible plus infected plus infectious populations:

$$S(t) + E(t) + I(t) = e^{-\mu t}\left(N - R(0) + \int_0^t e^{\mu \tau}((\omega + g)R(\tau) + (\mu - g_1)N)d\tau\right) = N - R(t) \qquad (44)$$

after using Assertion 1 so that one gets (38) from (40) as $t \to \infty$ which leads to

$$\lim_{t \to \infty} \int_0^t e^{-\mu(t-\tau)}(\omega + g)R(\tau)d\tau = \frac{(g_1 - \mu)N}{\mu} \qquad (45)$$

if $g_1 = \mu + \omega + g$ (implying that $S(\infty) = I(\infty) = R(\infty) = 0$ and $R(\infty) = N$). Finally, since $R(t) \leq N$ for all time and since $g_1 - g = \mu + \omega$ from (35), one gets from (36) for all time:

$$V(t) \geq 0 \Rightarrow I(t) \leq \frac{(g_1 - g)N}{\gamma} = \frac{(\mu + \omega)N}{\gamma} \leq \frac{g_1 N - gR(t)}{\gamma}$$

which is guaranteed for arbitrary initial conditions of (1)-(4) from Assertion 1 if $g_1 \geq \sigma$ and, in particular, if $\sigma = \mu + \omega$ and $g \geq 0$. $\qquad \square$

A further result for nonnegative vaccination being guaranteed for all time for the vaccination strategy of Theorem 3 is the following:

**Corollary 2**. The vaccination strategy of Theorem 3 is nonnegative for all time if

$$min\left(1, \frac{\sigma \beta}{(\mu + \sigma)(\mu + \gamma)}\right) \geq \frac{\mu + \omega}{\gamma}$$

*Proof*: Note $I(t) = \frac{\sigma E(t)}{\mu + \gamma}$ yields $\dot{I}(t) = 0$ from (3) so that a maximum or minimum of the infectious population is reached depending on the infected E(t). A potential maximum is reached for $I(t) = \frac{\sigma E(t)}{\mu + \gamma}$ with $E(t) \neq 0$. Thus, $I(t) \leq I_{max} := max(I(t) : t \geq 0) \leq \frac{\sigma E_{max}}{\mu + \gamma}$ where $E_{max} := \max(E(t): t \geq 0)$, is reached for $\dot{E}(t) = 0$ in (2) so that $\beta \frac{S_{max} I_{max}}{N} \geq (\mu + \sigma) E_{max}$, where $S_{max} := \max(S(t): t \geq 0)$. Combining the two relations yields the proof by using Assertion 1 since one has for all time:

$$I(t) \leq I_{max} \leq \frac{\sigma E_{max}}{\mu + \gamma} \leq \frac{\sigma \beta S_{max} I_{max}}{(\mu + \sigma)(\mu + \gamma) N} \leq \frac{\sigma \beta N}{(\mu + \sigma)(\mu + \gamma)}$$

$$\Rightarrow I(t) \leq I_{max} \leq \min\left(1, \frac{\sigma \beta}{(\mu + \sigma)(\mu + \gamma)}\right) N \qquad \square$$

**Remark 5**. An important problem to validate the SEIR model (1)-(4) under vaccination for practical application is the design of a vaccination strategy such that the obtained model is a positive system, as the real problem it describes is, in the sense that none of the populations S(t), E(t), I(t) and R(t) becomes negative at any time. It is proven in Appendix B (Theorem B.1) that if the vaccination strategy is constrained to the real interval $[0,1]$ for all time then none of those populations is negative at any time provided that all of them are nonnegative at t=0. Conditions to maintain its value under the positive unit are discussed in the following. $\square$

**Theorem 4**. Consider the vaccination strategy (35)-(36) of Theorem 3 subject to:

$$g_1 = \mu + \omega + g \;\; ; \;\; g < 0; \;\; \mu \geq |g| - \omega + \max(\gamma, |g|) \geq \max(\gamma, |g|) \qquad (46)$$

Thus, if all the partial populations of susceptible, infected, infectious and immune have nonnegative initial conditions then the vaccination function fulfills $V: \mathbf{R}_{0+} \to [0,1]$ provided that $(|g| - \omega)$ is sufficiently large. Also, all the values taken by any of those populations in the mathematical SEIR model are nonnegative for all time.

*Proof*: Since $\mu + \omega \geq |g| + \max(\gamma, |g|)$ then $g_1 = \mu + \omega + g \geq \max(\gamma, |g|) > 0$. Thus,

$$g_1 N - g R - \gamma I = g_1 R + g_1 I + g_1 (S + E) - g R - \gamma I = (g_1 - g) R + (g_1 - \gamma) I + g_1 (S + E)$$

$$\geq \min(g_1 - g, g_1 - \gamma, g_1)(R + I + S + E) = \min(g_1 + |g|, g_1 - \gamma, g_1) N \geq 0 \qquad (47)$$

and $V(t) \geq 0$ for all time from (36). The third constraint of (46) implies $\omega \leq |g|$ so that $V(t) > 1$ for some time t if

$$g_1 N + |g| R(t) - \gamma I(t) > \mu N$$

that is if

$$-(|g| - \omega) N + |g| R(t) - \gamma I(t) > 0$$

what is impossible for sufficiently large $(|g| - \omega)$. Since $V: \mathbf{R}_{0+} \to [0,1]$ then all the populations of the SEIR model are guaranteed to be nonnegative for all time from Theorem B.1. $\square$

**5 Feedback linearization in vaccination design**

The vaccination control laws proposed and analyzed in Theorems 2 and 3 can be regarded as special cases of a general design methodology called *feedback linearization*, [17, 18]. Feedback linearization is a general design methodology which has been successfully used in many non-linear control problems, [17-20]. The objective of this section is to frame the control laws introduced in the previous sections applied to the nonlinear epidemic model (1)-(4) into the feedback linearization formalism in order to present the rationale behind the proposed vaccination control laws, discuss some technical details concerning them and present the complete technique to be used in different epidemic models or vaccination control design. Thus, the general formalism presented, for instance, in [17, 18] is applied to the nonlinear system (1)-(4). The method requires us to follow a number of steps, [18] as follows: i) Initially, the relative degree of the system has to be obtained ii) Then, the nonlinear system (1)-(4) is to be re-written in a *normal canonical form.* iii) Next, the zero dynamics of the system are needed to be calculated and proved to be stable. This is a technical requirement on the system in order to guarantee the well-posededness of the control law. iv) Once the zero dynamics are proved to be stable, the design of the control law is direct from the canonical normal form. These steps are developed with detail in the following sections.

*5.1 Relative degree and normal canonical form*

The starting point is the epidemic model given by equations (1)-(4), re-written for convenience as:

$$\dot{x}(t) = f(x) + h(x)V(t) \tag{48}$$

with $x(t)^T = [S(t) \quad E(t) \quad I(t) \quad R(t)]$, $\beta' = \beta/N$ (which is a number since $N$ is a constant) and

$$f(x) = \begin{bmatrix} -\mu S(t) + \omega R(t) - \beta' S(t)I(t) + \mu N \\ \beta' S(t)I(t) - (\mu + \sigma)E(t) \\ -(\mu + \gamma)I(t) + \sigma E(t) \\ -(\mu + \omega)R(t) + \gamma I(t) \end{bmatrix}; \quad h(x) = \mu N \begin{bmatrix} -1 \\ 0 \\ 0 \\ 1 \end{bmatrix} \tag{49}$$

Along with state-space system (48)-(49), we also need to consider an output $y(t)$. The choice of different outputs for the system (48)-(49) leads to different control laws. This is an important fact for this method: a collection of control laws can be generated within this frame by just selecting different outputs $y(t)$. To illustrate the nature of the method, the output is selected as:

$$y(t) = R(t) \tag{50}$$

The next step is to calculate the relative degree of system (48)-(50). The relative degree can be defined as the number of times the output has to be derived until the input $V(t)$ appears in the derivative.

Thus, we can now derive the output equation (50):

$$\dot{y}(t) = \dot{R}(t) = -(\mu + \omega)R(t) + \gamma I(t) + \mu N V \tag{51}$$

It can be appreciated in Eq. (51) that the input $V(t)$ appears in the first derivative of $y(t)$. Thus, the system possess relative degree unity since the first derivative is enough to obtain the input. Furthermore, since $\mu N \neq 0$, the relative degree of system (48)-(50) is well defined in the complete state-space. The relative degree knowledge allows us to obtain a normal form for the original system. The normal form is a change of coordinates in the state-space that will permit the design of the feedback control law. The basic objective of the coordinates transformation is to obtain a nonlinear system with the input ($V(t)$) appearing in just one equation. According to [18], the first coordinate of the transformation is defined directly as the output:

$$z_1(t) = y(t) = R(t) \tag{52}$$

while the remaining variables, $\{z_i(t)\}_{i=2}^{4}$, are selected to satisfy the condition:

$$\left(\frac{\partial z_i(t)}{\partial x}\right)^T h(x) = \mu N \frac{\partial z_i(t)}{\partial S} - \mu N \frac{\partial z_i(t)}{\partial R} = 0 \tag{53}$$

for $i = 2,3,4$. Equation (53) becomes:

$$\frac{\partial z_i(t)}{\partial S} = \frac{\partial z_i(t)}{\partial R} \tag{54}$$

whose solution is given by:

$$z_i(t) = \lambda(E(t), I(t))\,(S(t) + R(t)); \qquad i = 2,3,4 \tag{55}$$

where $\lambda(E,I)$ is an arbitrary differentiable function of $(E,I)$. For the sake of simplicity, let's take $\lambda(E,I) = 1$, being the variable $z_2$ defined as:

$$z_2(t) = S(t) + R(t) \tag{56}$$

The remaining variables should be selected to satisfy Eq. (54) while being linearly-independent with (56). Fortunately, the seek of linearly independent solutions to (54) is not necessary in this case since equations (2) for $\dot{E}$ and (3) for $\dot{I}$ are not already directly dependent on the input, which is the objective we wanted to fulfill. Hence, it is made:

$$z_3 = E; \quad z_4 = I \tag{57}$$

The coordinates transformation (52)-(56)-(57) defines a global diffeomorphism in the state-space since

$$\frac{\partial(S,E,I,R)}{\partial(z_1,z_2,z_3,z_4)} = -1 \neq 0$$

and, therefore, the transformation is well-defined. This coordinates transformation converts the original system (49)-(50) into the system in normal form:

$$\dot{z}_1(t) = -(\mu+\omega)z_1(t)(t) + \gamma z_4(t) + \mu N V(t) \tag{58}$$

$$\dot{z}_2(t) = -\mu z_2(t)(t) + \gamma z_4(t) - \beta'(z_2 - z_1)z_4(t) + \mu N \tag{59}$$

$$\dot{z}_3(t) = \beta'(z_2 - z_1)z_4(t) - (\mu+\sigma)z_3 \tag{60}$$

$$\dot{z}_4(t) = -(\mu+\omega)z_4 + \sigma z_3 \tag{61}$$

$$y(t) = z_1(t) \tag{62}$$

Notice that the input only appears in the first equation of the system (58)-(61). This fact allows the design of the vaccination control. The next step is to analyze the zero dynamics of system (58)-(62).

### 5.2 Zero-dynamics of the normal system

This section analyzes the zero dynamics of system (58)-(62) which corresponds to the second step in the general process. The stability of the zero dynamics is crucial to ensure the applicability and stability of the vaccination strategy. The zero-dynamics can be regarded as the nonlinear counterpart of the zeros of a linear system and they are defined based on the output zeroing problem. This problem consist in finding an input signal which renders $z_1(t) = \dot{z}_1(t) = 0$ for all $t \geq 0$. Thus, from Eq. (58), such an input is defined by $V(t) = \frac{\gamma}{\mu N} z_4(t)$ which converts the first equation into the trivial one $0 = 0$. As far as the rest of variables, $z_2(t), z_3(t), z_4(t)$ are concerned, the system of equations (59)-(61) becomes:

$$\dot{z}_2(t) = -\mu z_2(t)(t) + \gamma z_4(t) - \beta' z_2(t) z_4(t) + \mu N \tag{63}$$

$$\dot{z}_3(t) = \beta' z_2(t) z_4(t) - (\mu + \sigma) z_3(t) \tag{64}$$

$$\dot{z}_4(t) = -(\mu + \omega) z_4(t) + \sigma z_3(t) \tag{65}$$

The set of equation (63)-(65) is said to be the *zero dynamics* of the nonlinear epidemic system (59)-(61). The following result holds:

**Lemma 1**. The zero dynamics of system (59)-(61) are stable, and thus, all variables $z_2(t), z_3(t), z_4(t)$ are bounded for all time.

*Proof*: The zero dynamics are defined by equations (63)-(65). Thus, summing up both sides of these equations:

$$\dot{z}_2(t) + \dot{z}_3(t) + \dot{z}_4(t) = \frac{d(z_2(t) + z_3(t) + z_4(t))}{dt} = 0 \tag{66}$$

implying that $z_2(t) + z_3(t) + z_4(t) = C$ (constant). In addition, it can be directly proved using similar arguments as those employed in Theorem B.1 that $z_i(t) \geq 0$ for $i = 2,3,4$ and all $t \geq 0$ (i.e. the zero-dynamics are positive). Consequently, $0 \leq z_i(t) \leq C$ for $i = 2,3,4$ and all $t \geq 0$ and the Lemma 1 is proved. □

In this way, the technical condition to guarantee the employment of feedback linearization control laws is satisfied and we are ready to design the vaccination strategy.

*5.3 Feedback control design*

The feedback control law *V(t)* is now designed by taking Eq. (58) and designing *V(t)* to cancel the dynamics of the right-hand terms of the equations in the form $V(t) = \frac{1}{\mu N}\left((\mu + \omega) z_1(t)(t) - \gamma z_4(t) + \eta(t)\right)$ which substituted into (58) yields $\dot{z}_1(t) = \eta(t)$ since all the terms appearing in Eq. (58) disappear due to the feedback control law. This cancellation plays the same role as the pole-zero cancellation in linear-systems which requires the stability of the zeros. The nonlinear counterpart of the linear zeros is the zero-dynamics whose stability has been verified to be stable in subsection 5.2 and thus, the feedback control makes the nonlinear system have no stability problems. Now, the signal $\eta(t)$ is used to govern the dynamics of $\dot{z}_1(t)$. A possible selection is $\eta(t) = -g' z_1 + g_1 N$ with $g', g_1 \geq 0$ aimed at making the immune match the total population, N. Thus, the complete control law (undoing the change of coordinates) becomes:

$$V(t) = \frac{1}{\mu N}\left((\mu + \omega) z_1(t)(t) - \gamma z_4(t) - g' z_1 + g_1 N\right) = \frac{1}{\mu N}\left((\mu + \omega - g') R(t)(t) - \gamma I(t) + g_1 N\right)$$
$$= \frac{1}{\mu N}\left(g R(t)(t) - \gamma I(t) + g_1 N\right) \tag{67}$$

which is exactly the vaccination law (36). Thus, the feedback linearization method is the analytical frame containing the presented control laws. A different choice for the output *y(t)* would lead to different vaccination strategies. For instance, the choice $y(t) = S(t)$ would lead to the vaccination law (26), but any other choice for the output could lead to an admissible vaccination policy. Thus, the general frame can generate a family of vaccination policies by selecting different choices for the output. The complete analysis of control law (67) has been performed in Section 4. It is worthwhile to note that despite being a well-known control design method for nonlinear systems, its application in epidemics has been rather limited. Simulation results have been performed for perfectly modeled casesand also under parametrical uncertainties in order to analyze the robustness of the proposal. The obtained results seems to be promising.

**Acknowledgements**

The authors thank to the Spanish Ministry of Education by its support of this work through Grant DPI2009-07197 and to the Basque Government by its support through Grants IT378-10 and SAIOTEK SPE07UN04.

**Appendix A** . *Equilibrium points and their stability properties*

*A.1) Equilibrium points of the uncontrolled system*

Assume for discussion simplicity of the equilibrium points that $\sigma = \gamma$ in the SEIR model (1)-(4). The equilibrium points $x^* = (S^*, E^*, I^*, R^*)^T$ of (1)-(4) under identically zero vaccination strategy satisfy the set of constraints:

$$\mu S^* - \omega R^* + \beta \frac{S^* I^*}{N} = \mu N \tag{A.1}$$

$$\beta \frac{S^* I^*}{N} = (\mu + \sigma) E^* \tag{A.2}$$

$$(\mu + \sigma) I^* = \sigma E^* \tag{A.3}$$

$$(\mu + \omega) R^* = \sigma I^* \tag{A.4}$$

An equilibrium point is $x_1^* = (N, 0, 0, 0)^T$. Another one is calculated as follows:

The combination of (A.2) –(A.3) yields:

$$\beta \frac{S^* I^*}{N} = \frac{(\mu + \sigma)^2 I^*}{\sigma} \Rightarrow S^* = \frac{(\mu + \sigma)^2}{\sigma \beta} N \text{ if } I^* \neq 0 \tag{A.5}$$

This constraint can be explored to obtain a new equilibrium point if $\frac{(\mu + \sigma)^2}{\sigma \beta} \leq 1$ guaranteeing the necessary model constraint $S^* \leq N$ from Assertion 1. If $\frac{(\mu + \sigma)^2}{\sigma \beta} \geq 1$ then the only equilibrium point is $x_1^*$ since $\frac{(\mu + \sigma)^2}{\sigma \beta} = 1 \Rightarrow S^* = N$. Thus, assume that $\frac{(\mu + \sigma)^2}{\sigma \beta} < 1$ and (A.5) holds. Then one gets from (A1) and (A.4) that:

$$\left[ \frac{(\mu + \sigma)^2 (\mu + \omega)}{\sigma^2} - \omega \right] R^* = \mu N \left( 1 - \frac{(\mu + \sigma)^2}{\sigma \beta} \right) \tag{A.6}$$

provided that $\mu > 0$ implying that $\frac{(\mu+\sigma)^2(\mu+\omega)}{\sigma^2} > \omega$. If $\mu = 0$ implying that $\frac{(\mu+\sigma)^2(\mu+\omega)}{\sigma^2} = \omega$ then the equilibrium points are $x_1^*$ and $x_{20}^* = \left( \frac{\sigma N}{\beta}, \frac{(\beta-\sigma)\omega N}{\beta(2\omega+\sigma)}, \frac{(\beta-\sigma)\omega N}{\beta(2\omega+\sigma)}, \frac{(\beta-\sigma)\sigma N}{\beta(2\omega+\sigma)} \right)^T$. Eq. (A.6) is equivalent to

$$R^* = \frac{\sigma\left(\sigma\beta - (\mu+\sigma)^2\right)}{\beta\left((\mu+\sigma)^2 + \omega(\mu+2\sigma)\right)} N \tag{A.7}$$

for $\mu > 0$ which implies $R^* \geq 0$ if $\frac{(\mu+\sigma)^2}{\sigma\beta} \leq 1$. From (A.4) and (A.7), one gets if $\sigma \neq 0$ that

$$I^* = \frac{\mu+\omega}{\sigma} R^* = \frac{(\mu+\omega)\left(\sigma\beta - (\mu+\sigma)^2\right)}{\beta\left((\mu+\sigma)^2 + \omega(\mu+2\sigma)\right)} N \tag{A.8}$$

Now, combining (A.3) and (A.8) yields:

$$E^* = \left(1 + \frac{\mu}{\sigma}\right) I^* = \frac{(\mu+\omega)(\mu+\sigma)}{\sigma^2} R^* = \frac{(\mu+\omega)(\mu+\sigma)\left(\sigma\beta - (\mu+\sigma)^2\right)}{\sigma\beta\left((\mu+\sigma)^2 + \omega(\mu+2\sigma)\right)} N \tag{A.9}$$

Thus,

$$x_2^* = \left( \frac{(\mu+\sigma)^2}{\sigma\beta} N, \frac{(\mu+\omega)(\mu+\sigma)\left(\sigma\beta-(\mu+\sigma)^2\right)}{\sigma\beta\left((\mu+\sigma)^2+\omega(\mu+2\sigma)\right)} N, \frac{(\mu+\omega)\left(\sigma\beta-(\mu+\sigma)^2\right)}{\beta\left((\mu+\sigma)^2+\omega(\mu+2\sigma)\right)} N, \frac{\sigma\left(\sigma\beta-(\mu+\sigma)^2\right)}{\beta\left((\mu+\sigma)^2+\omega(\mu+2\sigma)\right)} N \right)^T$$

(A.10)

is an equilibrium point of (1)–(4) provided that none of its components exceeds N and $\frac{(\mu+\sigma)^2}{\sigma\beta} < 1$ holds, that is if

$$\frac{\sigma\beta - (\mu+\sigma)^2}{\beta\left((\mu+\sigma)^2 + \omega(\mu+2\sigma)\right)} \max\left(\sigma, \left(1+\frac{\mu}{\sigma}\right)(\mu+\omega)\right) \leq 1 \tag{A.11}$$

□

*Remarks in the proof.* Note that $\frac{(\mu+\sigma)^2}{\sigma\beta} = 1$ then an equilibrium point $x_{21}^* = x_1^*$ exists as a particular case of (A.10). Also, $(N, 0, 0, 0)^T$ is an equilibrium point if $\mu = 0$ and $\beta = \sigma$ and $x_{20}^*$ is an equilibrium point if $\mu = 0$ and $\sigma < \beta$ both being particular cases of (A.10). Note also that both constraints of (A.11) are guaranteed in particular for sufficiently small positive parameters $\sigma = \gamma$ and μ after fixing $\omega, \beta$. □

*A.2) Stability of the linearized model about the equilibrium points*

The linearized model (1)-(4) about its equilibrium points is:

$$\begin{bmatrix} \Delta \dot{S}(t) \\ \Delta \dot{E}(t) \\ \Delta \dot{I}(t) \\ \Delta \dot{R}(t) \end{bmatrix} = \begin{bmatrix} -\mu - \beta \frac{I^*}{N} & 0 & -\beta \frac{S^*}{N} & \omega \\ \beta \frac{I^*}{N} & -(\mu+\sigma) & \beta \frac{S^*}{N} & 0 \\ 0 & \sigma & -(\mu+\sigma) & 0 \\ 0 & 0 & \sigma & -(\mu+\omega) \end{bmatrix} \begin{bmatrix} \Delta S(t) \\ \Delta E(t) \\ \Delta I(t) \\ \Delta R(t) \end{bmatrix} \quad (A.12)$$

At the equilibrium point $x_1^*$, the linearized system (A.12) becomes:

$$\begin{bmatrix} \Delta \dot{S}(t) \\ \Delta \dot{E}(t) \\ \Delta \dot{I}(t) \\ \Delta \dot{R}(t) \end{bmatrix} = \begin{bmatrix} -\mu & 0 & -\beta & \omega \\ 0 & -(\mu+\sigma) & \beta & 0 \\ 0 & \sigma & -(\mu+\sigma) & 0 \\ 0 & 0 & \sigma & -(\mu+\omega) \end{bmatrix} \begin{bmatrix} \Delta S(t) \\ \Delta E(t) \\ \Delta I(t) \\ \Delta R(t) \end{bmatrix} \quad (A.13)$$

whose characteristic equation is

$$p(s) = (s+\mu)\left[(s+\mu+\omega)(s+\mu+\sigma)^2 + \beta \, det\left(\begin{bmatrix} -\sigma & 0 \\ 0 & s+\mu+\omega \end{bmatrix}\right)\right]$$

$$= (s+\mu)(s+\mu+\omega)\left((s+\mu+\sigma)^2 - \sigma\beta\right) = 0$$

The characteristic zeros are $-\mu, -(\mu+\omega)$ and $-(\mu+\sigma\pm\sqrt{\sigma\beta})$. As a result, the equilibrium point $x_1^*$ of (A.13) is locally asymptotically Lyapunov stable if $\mu > 0, \omega > -\mu$ and $0 \leq \beta < \frac{(\mu+\sigma)^2}{\sigma}$. For the equilibrium point $x_2^*$, the linearized system (A.12) has a characteristic equation

$$p(s) = \left(s+\mu+\beta\frac{I^*}{N}\right) det\left(\begin{bmatrix} s+\mu+\sigma & -\beta\frac{S^*}{N} & 0 \\ -\sigma & s+\mu+\sigma & 0 \\ 0 & -\sigma & s+\mu+\omega \end{bmatrix}\right) + \omega \, det\left(\begin{bmatrix} -\beta\frac{I^*}{N} & s+\mu+\sigma & -\beta\frac{S^*}{N} \\ 0 & -\sigma & s+\mu+\sigma \\ 0 & 0 & -\sigma \end{bmatrix}\right)$$

$$+ \beta^2 \sigma(s+\mu+\omega)\frac{S^* I^*}{N^2}$$

$$= \left(s+\mu+\beta\frac{I^*}{N}\right)(s+\mu+\omega)\left((s+\mu+\sigma)^2 - \beta\sigma\frac{S^*}{N}\right) + \beta^2\sigma(s+\mu+\omega)\frac{S^* I^*}{N^2} - \beta\omega\sigma^2\frac{I^*}{N}$$

$$= p_0(s) + \beta \tilde{p}(s) = p_0(s)\left(1 + \frac{\beta \tilde{p}(s)}{p_0(s)}\right) = 0 \quad (A.14)$$

where

$$p_0(s) := (s+\mu)(s+\mu+\sigma)^2(s+\mu+\omega) \; ; \; \tilde{p}(s) := \frac{I^*}{N}\left((s+\mu+\sigma)^2(s+\mu+\omega) - \omega\sigma^2\right) - \sigma\frac{S^*}{N}(s+\mu)(s+\mu+\omega)$$

(A.15)

evaluated at (A.10). From the root locus technique in (A.14), the zeros of $p(s)$ converge to those of $p_0(s)$, namely, $s = -\mu$, $s = -(\mu+\sigma)$ (double), $s = -(\mu+\omega)$ as $\beta \to 0$. As a result, the eigenvalues of the linearized system (A.12) about $x_2^*$ are all stable from the continuity of the root locus for $|\beta|$ not exceeding some sufficiently small threshold value for any given values of the remaining parameters of (1)-(4). Equivalently, that property holds if

$$\left\| \frac{\tilde{p}(s)}{p_0(s)} \right\|_\infty := \max_{\omega \in \mathbf{R}_{0+}} \left| \frac{\tilde{p}(i\omega)}{p_0(i\omega)} \right| < \frac{1}{\beta}$$

where $\|.\|_\infty$ is the $\mathbf{RH}_\infty$-norm of strictly stable transfer functions, $i = \sqrt{-1}$ is the complex unit. This follows since $p_0(s)$ being a Hurwitz polynomial implies that $p(s)$ is Hurwitz if $|\beta \tilde{p}(i\omega)| < |p_0(i\omega)|$; $\forall \omega \in \mathbf{R}_{0+}$ from Rouché theorem of number of zeros within a closed set applied to the complex half-plane $Re\, s < 0$. Note that the global Lyapunov stability is automatically guaranteed for the SEIR model (1)-(4) since the total population is assumed to be constant for all time. The summarized local stability result around the equilibrium points is as follows:

**Theorem A.1**. The vaccination-free SEIR model (1)-(4) is locally stable about $x_1^*$. It is also locally stable about $x_2^*$ if $\left\| \dfrac{\tilde{p}(s)}{p_0(s)} \right\|_\infty < \dfrac{1}{\beta}$. □

**Appendix B**. *Positive solutions of the SEIR model (1)-(4)*

The following result holds:

**Theorem B.1**. Assume the SEIR model (1)-(4) with $N = N(0) = S(0) + E(0) + I(0) + R(0) > 0$ and $min(S(0), E(0), I(0), R(0)) \geq 0$ under any vaccination strategy $V : \mathbf{R}_{0+} \to [0,1]$. Then, $min(S(t), E(t), I(t), R(t)) \geq 0$; $\forall t \in \mathbf{R}_{0+}$.

*Proof*: Let eventually exist finite time instants $t_S, t_E, t_I, t_R \in \mathbf{R}_{0+}$ with $t^* := min(t_S, t_E, t_I, t_R)$ be such that:

If $t^* = t_S$ then $S(t_S) = 0$, $min(S(t), E(t), I(t), R(t)) \geq 0$; $\forall t \in [0, t_S]$

If $t^* = t_E$ then $E(t_E) = 0$, $min(S(t), E(t), I(t), R(t)) \geq 0$; $\forall t \in [0, t_E]$

If $t^* = t_I$ then $I(t_I) = 0$, $min(S(t), E(t), I(t), R(t)) \geq 0$; $\forall t \in [0, t_I]$

If $t^* = t_R$ then $R(t_R) = 0$, $min(S(t), E(t), I(t), R(t)) \geq 0$; $\forall t \in [0, t_R]$

Note that either $t^*$ does not exist or it is the first eventual finite time instant where some of the partial populations of the SEIR model reaches a zero value and can be coincident with at most three of its arguments since the total population being N is incompatible with the four partial populations being zero. The remaining of the proof is split into four parts as follows:

a) Proceed by contradiction by assuming that there exists a finite $t^* = t_S \geq 0$ such that $S(t) \geq 0$; $\forall t \in [0, t_S)$, $S(t_S) = 0$ and $S(t_S^+) < 0$, meaning with abbreviate notation that $S(t) < 0$;

$\forall t \in (t_S + \varepsilon_1, t_S + \varepsilon_1 + \varepsilon_2)$, with $min(E(t), I(t), R(t)) \geq 0$; $\forall t \in [0, t_S]$. Thus, $\dot{S}(t_S) = \omega R(t_S) + \mu N(1 - V(t_S)) \geq 0$ from (1) since $V(t) \in [0,1]$; $\forall t \in \mathbf{R}_{0+}$. Since $S(t_S) = 0$ and $\dot{S}(t_S) \geq 0$ then $S(t_S^+) \geq 0$, meaning with abbreviate notation that $S(t) > 0$; $\forall t \in (t_S + \varepsilon_1, t_S + \varepsilon_1 + \varepsilon_2)$, since the solution of the SEIR model (1)-(4) is continuous for all time, contradicting the assumption $S(t_S^+) < 0$ so that such a time instant $t^* = t_S \geq 0$ does not exist.

b) Proceed by contradiction by assuming that there exists a finite $t^* = t_E \geq 0$ such that $E(t) \geq 0$; $\forall t \in [0, t_E)$, $E(t_E) = 0$ and $E(t_E^+) < 0$ with $min(S(t), I(t), R(t)) \geq 0$; $\forall t \in [0, t_E]$. Thus, $\dot{E}(t_E) = \frac{\beta S(t_E) I(t_E)}{N} \geq 0$ from (2); $\forall t \in \mathbf{R}_{0+}$. Since $E(t_E) = 0$ and $\dot{E}(t_E) \geq 0$ then $E(t_E^+) \geq 0$, since the solution of the SEIR model (1)-(4) is continuous for all time, contradicting the assumption $E(t_E^+) < 0$ so that such a $t^* = t_E \geq 0$ does not exist.

c) Proceed by contradiction by assuming that there exists a finite $t^* = t_I \geq 0$ such that $I(t) \geq 0$; $\forall t \in [0, t_I)$, $I(t_I) = 0$ and $I(t_I^+) < 0$ with $min(S(t), E(t), R(t)) \geq 0$; $\forall t \in [0, t_I]$. Thus, $\dot{I}(t_I) = \sigma E(t_I) \geq 0$ from (3); $\forall t \in \mathbf{R}_{0+}$. Since $I(t_I) = 0$ and $\dot{I}(t_I) \geq 0$ then $I(t_I^+) \geq 0$, since the solution of the SEIR model (1)-(4) is continuous for all time, contradicting the assumption $I(t_I^+) < 0$ so that such a $t^* = t_I \geq 0$ does not exist.

d) Proceed by contradiction by assuming that there exists a finite $t^* = t_R \geq 0$ such that $R(t) \geq 0$; $\forall t \in [0, t_R)$, $R(t_R) = 0$ and $R(t_R^+) < 0$ with $min(S(t), E(t), I(t)) \geq 0$; $\forall t \in [0, t_R]$. Thus, $\dot{R}(t_R) = \gamma I(t_R) + \mu N V(t_R) \geq 0$ from (4) since $V(t) \in [0,1]$; $\forall t \in \mathbf{R}_{0+}$. Since $R(t_R) = 0$ and $\dot{R}(t_R) \geq 0$ then $R(t_R^+) \geq 0$, since the solution of the SEIR model (1)-(4) is continuous for all time, contradicting the assumption $R(t_R^+) < 0$ so that such a time instant $t^* = t_R \geq 0$ does not exist.

If such a finite time instant $t^*$ does not exist then the above result follows directly. As a result, $min(S(0), E(0), I(0), R(0)) \geq 0 \Rightarrow min(S(t), E(t), I(t), R(t)) \geq 0$; $\forall t \in \mathbf{R}_{0+}$ since there is no time instant $t^*$ for which any of the four partial populations reaches a zero value with its first time-derivative being simultaneously negative at such a time instant. □